\documentclass[a4paper, 12pt] {article}
\usepackage[cp1251]{inputenc}
\usepackage[english]{babel}
\usepackage{amsfonts}
\usepackage{amsmath}
\usepackage{mathtext}
\usepackage{cite}
\newcommand{\de}{\delta}

\newcommand{\al}{\alpha}

\newcommand{\ee}{\varepsilon}

\begin{document}

\begin{center}
{\large\bf On inverse spectral problems for Sturm--Liouville differential operators on closed sets}\\[0.4cm]
{\bf S.A.\,Buterin, M.A.\,Kuznetsova and  V.A.\,Yurko} \\[0.4cm]
\end{center}

\thispagestyle{empty}

{\bf Abstract.} We study Sturm--Liouville operators on closed sets of a special structure, which are sometimes referred as time scales and
often appear in modelling various real processes. Depending on the set structure, such operators unify both differential and difference
operators. We obtain properties of their spectral characteristics and study inverse problems with respect to them. We prove that the
spectral characteristics uniquely determine the operator.

{\it Keywords:} differential operators; closed sets; time scales; inverse spectral problems.

{\it AMS Mathematics Subject Classification (2010):} 34A55 34B24 47E05\\

{\bf 1. Introduction}
\\

Differential operators on a closed set of the real line, which is sometimes called the time scale, often appear in modelling various real
processes (see \cite{[1],[2]}). Depending on the set structure, such operators unify derivatives and differences. We study inverse spectral
problems for Sturm--Liouville differential operators on a closed subset of the real line. Inverse spectral problems consist in recovering
operators from given spectral characteristics. For the classical Sturm--Liouville operators on an interval, inverse problems have been
studied fairly completely; the classical results can be found in \cite{[3], [4], [5], [6]}. However, nowadays there is no inverse problem
theory for differential operators defined on general closed sets. The only one exception is the work \cite{[7]} where an Ambarzumian-type
theorem is proved for Sturm--Liouville differential operators on time scales.

In this paper we consider a subclass of closed sets that are finite unions of non-overlapping segments. To determine the Sturm--Liouville
operator on such sets, we provide main notions of the time scale theory in Section 2. In Section~3, we study properties of spectral
characteristics including the Weyl function. In Section~4, we study the inverse problems of recovering the potential of the
Sturm--Liouville operator from the given Weyl function as well as the spectra of two boundary value problems with one common boundary
condition. The uniqueness theorems for these inverse problems are proved (Theorems~1 and~2).
\\

{\bf 2. The main notions}
\\

First of all, we introduce several notions of the time scale theory (see \cite{[1],[2]} for more details). Let $T$ be a closed subset of
${\mathbb R},$ which we refer to as the time scale. We define the so-called jump functions $\sigma$ and $\sigma_-$ on $T$ in the following
way:
$$
\sigma(x)=\left\{\begin{array}{cl}\inf \{s\in T:\; s>x\}, & x\ne \sup T,\\[2mm]
\sup T, & x=\sup T,
\end{array}\right.
\sigma_-(x)=\left\{\begin{array}{cl}\sup \{s\in T:\; s<x\}, & x\ne \inf T,\\[2mm]
\inf T, & x=\inf T.
\end{array}\right.
$$
A point $x\in T$ is called {\it left-dense}, {\it left-isolated}, {\it right-dense} and {\it right-isolated}, if $\sigma_-(x)=x,$
$\sigma_-(x)<x,$ $\sigma(x)=x$ and $\sigma(x)>x,$ respectively. If $\sigma_-(x)<x<\sigma(x),$ then $x$ is called {\it isolated}; if
$\sigma_-(x)=x=\sigma(x),$ then $x$ is called {\it dense}. Denote $T^0:=T\setminus\{\sup T\},$ if $\sup T$ is left-isolated, and $T^0:=T,$
otherwise. We also denote by $C(B)$ the class of continuous functions on the subset $B \subseteq T.$

A function $f$ on $T$ is called $\Delta$-{\it differentiable} at $t \in T^0,$ if for any $\ee>0$ there exists $\delta>0$ such that
$$
|f(\sigma(t))-f(s)-f^{\Delta}(t)(\sigma(t)-s)|\le\ee |\sigma(t)-s|
$$
for all $s\in (t-\de, t+\de)\cap T.$ The value $f^{\Delta}(t)$ is called the $\Delta$-{\it derivative} of the function $f$ at the
point~$t.$ The following proposition gives conditions of $\Delta$-differentiability at points of different types.

\medskip
{\bf Proposition 1. }{\it 1)  If $f(t)$ is $\Delta$-differentiable at $t,$ then $f(t)$ is continuous in $t.$

2) Let $t \in T$ be a right-isolated point.  Then $f$ is $\Delta$-differentiable at $t,$ if and only if $f$ is continuous in $t.$ In this
case we have
$$
f^{\Delta}(t)=\frac{f(\sigma(t))-f(t)}{\sigma(t)-t}.
$$
In particular, if $T=\{x=hk:\; k\in{\mathbb Z}\},$ then
$$
f^{\Delta}(x)=\frac{f(x+h)-f(x)}{h}.
$$

3) Let $t\in T$ be a right-dense point. Then $f$ is  $\Delta$-differentiable at $t,$ if and only if there exists the limit
$$
\lim_{s\to t,\; s \in T} \frac{f(t)-f(s)}{t-s}=:  f^{\Delta}(t).
$$

In particular, if $(t - \ee, t + \ee) \subset T$ for some $\ee > 0,$ then $f$ is $\Delta$-differentiable at $t,$ if and only if $f$ is
differentiable at $t.$ In this case the equality $f^{\Delta}(t)=f'(t)$ is fulfilled.}

\medskip
We also introduce derivatives of the higher order $n \ge 2.$ Let the $(n-1)$-th $\Delta$-derivative $f^{\Delta^{n-1}}$ of $f$ be defined on
$T^{0^{n-1}},$ where $a^n = \underbrace{a \ldots a}_{n}$ for any symbol $a.$ If $f^{\Delta^{n-1}},$ in turn, is $\Delta$-differentiable on
$T^{0^n}:= (T^{0^{n-1}})^0,$ then $f^{\Delta^n}:= (f^{\Delta^{n-1}})^{\Delta}$ is called the {\it $n$-th $\Delta$-derivative} of $f$ on
$T^{0^n}.$ For $n\ge1$ we also denote by $C^n(T)$ the class of functions $f$ for which there exists the $n$-th $\Delta$-derivative
$f^{\Delta^n}$ and $f^{\Delta^n} \in C(T^{0^n}).$

Consider the Sturm--Liouville equation on $T:$
\begin{equation}
\ell y:=-y^{\Delta\Delta}(x)+q(x)y(\sigma(x))=\lambda y(\sigma(x)), \quad x\in T^{0^2}, \label{1}
\end{equation}
where $\lambda$ is the spectral parameter and $q(x)\in C(T)$ is a complex-valued function. A function $y$ is called a solution of equation
(\ref{1}), if $y\in C^2(T^{0^2})$ and satisfies equation (\ref{1}). The statement and the study of inverse spectral problems essentially
depend on the structure of the time scale $T.$ It is necessary to choose and describe subclasses of time scales for which the inverse
problem theory can be constructed adequately. In this paper we consider the so-called $Y2$-structure, which has the form
$$
T=\bigcup_{k=1}^N \,[a_k, b_k],\quad N\ge 2,\quad a_k<b_k<a_{k+1}<b_{k+1}, \quad k=\overline{1,N-1}.
$$

For $j=0,1$ denote by $L_j$ the boundary value problem for equation (\ref{1}) on $T$ with the boundary conditions
$y^{\Delta^j}(a_1)=y(b_N)=0.$ Let $S(x,\lambda)$ and $C(x,\lambda)$ be solutions of equation~(\ref{1}) on $T$ satisfying the initial
conditions
$$
S(a_1,\lambda)=C^{\Delta}(a_1,\lambda)=0,\quad S^{\Delta}(a_1,\lambda)=C(a_1,\lambda)=1.
$$
One can show that for each fixed $x\in T$ the functions $S(x,\lambda)$ and $C(x,\lambda)$ are entire in $\lambda$ of order $1/2.$ We
introduce the entire functions
$$
\Delta_0(\lambda):=S(b_N,\lambda), \quad \Delta_1(\lambda) := C(b_N, \lambda),
$$
which are called the {\it characteristic functions} of the problems $L_0$ and $L_1,$ respectively. For $j=0,1$ eigenvalues
$\{\lambda_{nj}\}_{n\ge 1}$ of the problem $L_j$ coincide with zeros of $\Delta_j(\lambda).$

Further, let $\Phi(x,\lambda)$ be a solution of equation (\ref{1}) on $T$ satisfying the boundary conditions
\begin{equation}
\Phi^\Delta(a_1,\lambda)=1,\quad \Phi(b_N,\lambda)=0.                               \label{8}
\end{equation}
The function $M(\lambda):=\Phi(a_1,\lambda)$ is called the Weyl-type function or simply the Weyl function. It can also be shown that
\begin{equation}
\Phi(x,\lambda)=S(x,\lambda)+M(\lambda)C(x,\lambda),                              \label{9}
\end{equation}
\begin{equation}
M(\lambda)=-\frac{\Delta_0(\lambda)}{\Delta_1(\lambda)}.                              \label{10}
\end{equation}

The Weyl function $M(\lambda)$ as well as the pair of spectra $\{ \lambda_{nj}\}_{n \ge 1}$ are called spectral characteristics of the
operator $\ell.$ In the next section we study their properties and in Section~4 we prove that their specifications uniquely determine the
potential $q(x)$ on $T.$

According to Proposition~1, the $\Delta$-derivative of the function $y$ on the $Y2$--structure $T$ has the form
\begin{equation}
y^{\Delta}(b_k)=\frac{y(a_{k+1})-y(b_k)}{a_{k+1}-b_k},\;\;k=\overline{1,N-1}, \quad y^{\Delta}(x)=y'(x),\;\;x\in [a_k,b_k),\;\;
k=\overline{1,N}.     \label{2}
\end{equation}
Since $y^\Delta \in C(T)$ for any solution $y$ of equation (\ref{1}), for $k=\overline{1,N}$ we obtain $y^{\Delta}(b_k)=y'(b_k),$ where
$y'(b_k)$ is the left derivative at $b_k,$ which, obviously, exists. Consequently,
\begin{equation}
y(a_{k+1})=y(b_k)+(a_{k+1}-b_k)y'(b_k),\quad k=\overline{1,N-1}.    \label{3}
\end{equation}
Equalities (\ref{1}), (\ref{2}) and (\ref{3}) yield $N$ classical Sturm--Liouville equations on intervals
\begin{equation}
-y''(x)+q(x)y(x)=\lambda y(x),\quad x\in (a_k,b_k),\quad k=\overline{1,N},  \label{4}
\end{equation}
along with the relations
$$
y^{\Delta\Delta}(b_k)=\frac{y'(a_{k+1})-y'(b_k)}{a_{k+1}-b_k}  =(q(b_k)-\lambda)y(a_{k+1}),\quad k=\overline{1,N-1}.
$$
According to (\ref{2}), the latter relations are equivalent to
\begin{equation}
y'(a_{k+1})=y'(b_k)+(a_{k+1}-b_{k})(q(b_{k})-\lambda)y(a_{k+1}), \quad k=\overline{1,N-1}.
\label{5}
\end{equation}
From (\ref{3}) and (\ref{5}) we obtain the following jump conditions on the function $y(x):$
\begin{equation}
\left. \begin{array}{c}
y(a_{k+1})=\al^{k}_{11}(\lambda)y(b_k)+\al^{k}_{12}(\lambda)y'(b_k), \\[3mm]
y'(a_{k+1})=\al^{k}_{21}(\lambda)y(b_k)+\al^{k}_{22}(\lambda)y'(b_k),
\end{array} \right\}  \quad    k=\overline{1,N-1},                                           \label{6}
\end{equation}
where
\begin{equation}
\left. \begin{array}{ll}
\al^{k}_{11}(\lambda)=1, & \al^{k}_{12}(\lambda)=a_{k+1}-b_k,\\[3mm]
\al^{k}_{21}(\lambda)=(a_{k+1}-b_k)(q(b_k)-\lambda), & \al^{k}_{22}(\lambda)=1+(a_{k+1}-b_k)^2(q(b_k)-\lambda).
\end{array} \right\}                                                \label{7}
\end{equation}
Thus, for $j=0,1$ the problem $L_j$ is equivalent to the boundary value problem $L_j^1$ for the set of equations (\ref{4}) subject to the
jump conditions (\ref{6}) and the boundary conditions $y^{(j)}(a_1)=y(b_N)=0.$ Moreover, the functions $S(x,\lambda),$ $C(x,\lambda)$ and
$\Phi(x,\lambda)$ are solutions of equations (\ref{4}), satisfy the conditions (\ref{6}) and
\begin{equation}    \label{7.1}
S(a_1,\lambda)=C'(a_1,\lambda)= \Phi(b_N,\lambda)=0,\quad S'(a_1,\lambda)=C(a_1,\lambda)= \Phi'(a_1,\lambda)=1,
\end{equation}
and determined uniquely.

The boundary value problems $L_j^1,\,j=0,1,$ remind the so-called discontinuous boundary value problems for Sturm--Liouville equations on
an interval. We note that various aspects of inverse problems for discontinuous Sturm--Liouville operators were studied in
\cite{Kru,Hal,Shep,Yur00-1,Yur00-2,FrYur-02,ShYur,Yang09,Wang13,Yang13,AmOz14,Wang15} and other works. However, the special dependence of
the coefficients in jump conditions (\ref{6}) on the potential as well as on the spectral parameter require a separate investigation.
\\

{\bf 3. Properties of the spectral characteristics}
\\

Without loss of generality in what follows we assume that $a_1=0.$ In order to use the classical theory of Sturm--Liouville operators on an
interval, it is convenient to extent the function $q(x)$ on the whole segment $[0,b_N]$ such that $q(x)\in C[0,b_N]$ and arbitrary in the
rest. Consider the Sturm--Liouville equation
\begin{equation}
-y''(x)+q(x)y(x)=\lambda y(x),\quad x\in [0,b_N].                     \label{11}
\end{equation}
Let $\lambda=\rho^2$ and $\mathrm{Im}\, \rho \ge 0.$ It is known (see, for example, \cite{[5]}) that there exists a fundamental system
$\{Y_1(x,\rho), Y_2(x,\rho)\}$ of solutions of equation (\ref{11}) on $[0,b_N],$ having the asymptotics
\begin{equation}\label{12}
Y_1^{(j)}(x,\rho)=(i\rho)^j\exp(i\rho x)[1],\quad Y_2^{(j)}(x,\rho)=(-i\rho)^j\exp(-i\rho x)[1],\quad j=0,1,
\end{equation}
uniformly for $x\in[0,b_N],$ where $[1]=1+O(\rho^{-1}),$ $\rho\to\infty.$ Moreover, for each fixed $x\in[0,b_N]$ the functions
$Y_\nu^{(j)}(x,\rho)$ are continuous for $\mathrm{Im}\,\rho\ge0$ and analytic for $\mathrm{Im}\,\rho>0.$

Expanding $S(x, \lambda)$ and $C(x, \lambda)$ with respect to the system  $\{Y_1(x,\rho), Y_2(x,\rho)\},$ we get
\begin{equation} \left. \begin{array}{cc}
S(x,\lambda)=X^0_{2k-1}(\rho)Y_1(x,\rho)+X^0_{2k}(\rho)Y_2(x,\rho), \\[3mm]
C(x,\lambda)=X^1_{2k-1}(\rho)Y_1(x,\rho)+X^1_{2k}(\rho)Y_2(x,\rho),
\end{array} \right\} \quad x\in[a_k,b_k],\quad k=\overline{1,N}, \label{24}
\end{equation}
where the vectors $X^j = \big(X^j_{k}(\rho)\big)_{k=1}^{2N},$ $j=0, 1,$ can be found from jump conditions (\ref{6}) and initial conditions
(\ref{7.1}). Indeed, substituting (\ref{24}) into (\ref{6}) and (\ref{7.1}) we obtain the systems
\begin{equation}\label{25}
B X^0 = (0,-1,0,\ldots,0)^T, \quad B X^1 = (-1,0,0,\ldots,0)^T,
\end{equation}
where $T$ is the transposition sign and $B$ is the $2N \times 2N$ matrix of the following structure:
$$
B = \begin{pmatrix}
r_{11} & s_{11} & 0 & 0 & 0 & 0 & 0 & \ldots & 0 & 0 & 0 &  0  & 0 \\
r_{12} & s_{12}  & 0 & 0 & 0 & 0 & 0 & \ldots & 0 & 0 & 0 & 0 & 0\\

p_{21} & q_{21} & r_{21} & s_{21} & 0 & 0 & 0 & \ldots & 0 & 0 & 0 &  0  & 0 \\
p_{22} & q_{22} & r_{22} & s_{22} & 0 & 0 & 0 & \ldots & 0 & 0 & 0 &  0  & 0 \\

0 & 0 & p_{31} & q_{31} & r_{31} & s_{31} &  0 & \ldots & 0 & 0 & 0 & 0  & 0 \\
0 & 0 & p_{32} & q_{32} & r_{32} & s_{32} &  0 & \ldots & 0 & 0 & 0 & 0  & 0  \\

\vdots & \vdots &\vdots &\vdots & \vdots & \vdots &  \vdots & \ddots & \vdots & \vdots & \vdots & \vdots & \vdots   \\

0 & 0 &0 &0 & 0 & 0 & 0 & \ldots & 0 & p_{N1} & q_{N1} & r_{N1} & s_{N1} \\
0 & 0 &0 &0 & 0 & 0 & 0 & \ldots & 0 & p_{N2} & q_{N2} & r_{N2} & s_{N2}
\end{pmatrix}
$$
with
\begin{equation}\label{7.1.1}\left.\begin{array}{c}
p_{k\nu} = \alpha^{k-1}_{\nu1}(\lambda)Y_1(b_{k-1}, \rho) + \alpha^{k-1}_{\nu2}(\lambda)Y'_1(b_{k-1}, \rho), \quad k=\overline{2, N},\\[3mm]
q_{k\nu} =\alpha^{k-1}_{\nu1}(\lambda)Y_2(b_{k-1}, \rho) + \alpha^{k-1}_{\nu2}(\lambda)Y'_2(b_{k-1}, \rho), \quad k=\overline{2, N},\\[3mm]
r_{k\nu} = -Y_1^{(\nu-1)}(a_k, \rho), \quad s_{k\nu} = -Y_2^{(\nu-1)}(a_k, \rho), \quad k=\overline{1, N},
\end{array}\right\}
\end{equation}
where $\nu=1,2.$ The other elements of $B$ equal to zero.

Due to the structure of the matrix $B,$ its determinant can be calculated by the formula
\begin{equation}\label{7.1.2}
\det B =\prod_{k=1}^N (r_{k1} s_{k2} - r_{k2} s_{k1}) = \left(Y'_2(0, \rho)Y_1(0, \rho) - Y'_1(0, \rho)Y_2(0, \rho)\right)^N =
(-2i\rho)^N[1].
\end{equation}
Using Cramer's formulae it is easy to show that
\begin{equation}\label{7.2}
\Delta_j(\lambda) = (-1)^{j+1}\frac{\det D_j}{\det B}, \quad j=0,1,
\end{equation}
where
$$
D_j = \begin{pmatrix}
r_{1,j+1} & s_{1,j+1}  & 0 & 0 & 0 & 0 & 0 & \ldots & 0 & 0 & 0 & 0 & 0\\

p_{21} & q_{21} & r_{21} & s_{21} & 0 & 0 & 0 & \ldots & 0 & 0 & 0 &  0  & 0 \\
p_{22} & q_{22} & r_{22} & s_{22} & 0 & 0 & 0 & \ldots & 0 & 0 & 0 &  0  & 0 \\

0 & 0 & p_{31} & q_{31} & r_{31} & s_{31} &  0 & \ldots & 0 & 0 & 0 & 0  & 0 \\
0 & 0 & p_{32} & q_{32} & r_{32} & s_{32} &  0 & \ldots & 0 & 0 & 0 & 0  & 0  \\

\vdots & \vdots &\vdots &\vdots & \vdots & \vdots &  \vdots & \ddots & \vdots & \vdots & \vdots & \vdots & \vdots   \\

0 & 0 &0 &0 & 0 & 0 & 0 & \ldots & 0 & p_{N1} & q_{N1} & r_{N1} & s_{N1} \\
0 & 0 &0 &0 & 0 & 0 & 0 & \ldots & 0 & p_{N2} & q_{N2} & r_{N2} & s_{N2} \\
0 & 0 & 0 & 0 & 0 & 0 & 0 & \ldots & 0 & 0 & 0 &  Y_1(b_N,\rho)  & Y_2(b_N,\rho)
\end{pmatrix}.
$$

For $k = \overline{1, N}$ we denote by ${\cal D}_j^{k}$ the minor of the matrix $D_j,$ consisting of intersections of the columns
$2(N-k)+1,2(N-k)+2, \ldots, 2N$ with the following rows:
$$
2 (N-k) + 1,\; 2(N-k)+2,\; 2(N-k) + 3\;\, \ldots,\; 2N, \quad \text{if } j=1 \text{ or } k=N,
$$
$$
2 (N-k),\; 2(N-k)+2,\; 2(N-k) + 3,\; \ldots,\; 2N, \quad \text{if } j=0 \text{ and } k<N.
$$
The value ${\cal D}_j^{k}$ is the determinant of a $2k\times 2k$ matrix. In particular, ${\cal D}_j^{N} = \det D_j,$ $j=0,1.$

\medskip
{\bf Lemma 1. }{\it For $\rho\to\infty$ and $Im\rho\ge0$ the following asymptotics holds:
\begin{equation}\label{27}
{\cal D}_j^k =\left\{\begin{array}{l} \!\!\displaystyle -2(i \rho)^j  \Big(f^1_{1-j}(\rho) + O\Big( \frac1\rho\exp(\mathrm{Im}
\rho\gamma_1)\Big)\Big), \quad k = 1,\\[4mm]
\!\!\displaystyle (-2)^k (i \rho)^{k+j-1} \prod_{l=1}^{k-1}\alpha^{N-l}_{22}(\lambda)\Big( f_j^k(\rho)f_1^1(\rho)
\prod_{l=2}^{k-1}f_0^l(\rho) + O\Big( \frac1\rho\exp(\mathrm{Im} \rho\gamma_k)\Big)\Big), \; k > 1,
\end{array}\right.
\end{equation}
where $f^k_0(\rho) =\cos \rho (b_{N-k+1} - a_{N-k+1}),$ $f^k_1(\rho) = -i \sin \rho (b_{N-k+1} - a_{N-k+1})$ and
$$
\gamma_k = \sum_{l=N-k+1}^N (b_l - a_l), \quad k = \overline{1, N}, \quad  \prod_{l=2}^1f_0^l(\rho):=1.
$$}

{\it Proof.} According the definition of ${\cal D}_j^k$ and (\ref{12}), we have
$$
{\cal D}^1_j= \begin{vmatrix}
-Y_1^{(j)}(a_N, \rho) & -Y_2^{(j)}(a_N, \rho) \\
Y_1(b_N, \rho) & Y_2(b_N, \rho)
\end{vmatrix} = (-i\rho)^j\exp(i\rho(b_N-a_N))[1] -(i\rho)^j\exp(-i\rho(b_N-a_N))[1]
$$
$$
\qquad\qquad\qquad =(-i\rho)^j\exp(i\rho(b_N-a_N)) -(i\rho)^j\exp(-i\rho(b_N-a_N)) +O(\rho^{j-1}\exp(\mathrm{Im}\rho\gamma_1)),
$$
which coincides with (\ref{27}) for $k=1.$

Assume by induction that formulae (\ref{27}) are fulfilled for $k=n < N.$ Then expanding ${\cal D}_j^{n+1}$ along its first row, we obtain
\begin{equation} \label{minors}
{\cal D}_j^{n+1} = r_{N-n, j+1} \big(q_{N-n+1,1} {\cal D}_1^{n} - q_{N-n+1,2} {\cal D}_0^{n}\big) - s_{N-n,j+1}\big(p_{N-n+1,1} {\cal
D}_1^{n} - p_{N-n+1,2} {\cal D}_0^{n}\big).
\end{equation}
By virtue of (\ref{7}) and (\ref{12}), for $\nu=1,2$ we have
\begin{equation}\label{alpha}
\frac{\alpha^{N-n}_{1\nu}(\lambda)}{\alpha^{N-n}_{2\nu}(\lambda)} =O\Big(\frac1{\rho^2}\Big), \quad \rho\to\infty,
\end{equation}
and
\begin{equation} \label{pq}
p_{N-n+1,\nu} =i\rho\alpha^{N-n}_{\nu2}(\lambda)\exp(i \rho b_{N-n})[1], \quad q_{N-n+1,\nu} =-i\rho\alpha^{N-n}_{\nu2}(\lambda)\exp(-i
\rho b_{N-n})[1],
\end{equation}
$$
r_{N-n,j+1} = -(i \rho)^{j}\exp(i \rho a_{N-n})[1], \quad s_{N-n,j+1} = -(-i \rho)^j\exp(-i \rho a_{N-n})[1].
$$
According to
(\ref{27}) and (\ref{minors}), these asymptotic formulae yield
\begin{multline*}
{\cal D}_j^{n+1} = (i \rho)^{j+1} \alpha^{N-n}_{22}(\lambda) \Big((-1)^{j+1}\exp(i\rho(b_{N-n} -a_{N-n})){\cal D}_0^n -\exp(-i\rho(b_{N-n}
-a_{N-n})) {\cal D}_0^n \\ +O(\rho^{3n-4}\exp(\mathrm{Im}\rho\gamma_{n+1}))\Big),
\end{multline*}
which gives (\ref{27}) for $k=n+1.$ $\hfill\Box$

\medskip
By virtue of (\ref{7.1.2}), (\ref{7.2}) and (\ref{27}) for $k=N,$ we get
\begin{equation}\label{Delta}
\Delta_j(\lambda) = (-i\rho)^{j-1} \prod_{l=1}^{N-1} \alpha^l_{22}(\lambda)\Big(f_j^N(\rho)f_1^1(\rho)\prod^{N-1}_{l=2}f_0^l(\rho)
+O\Big(\frac1\rho\exp(\mathrm{Im \rho}\gamma_N)\Big)\Big), \quad \rho\to\infty.
\end{equation}
Using asymptotics (\ref{Delta}), in the standard way involving Rouch\'e's theorem one can establish that for $j=0,1$ the characteristic
function $\Delta_j(\lambda)$ possesses countably many zeros $\{\lambda_{nj}\}_{n\ge1}.$ Moreover, the following asymptotics holds:
$$
\{\lambda_{nj}\}_{n\ge N+j} =\bigcup_{\nu=1}^N \{\lambda_{k\nu j}\}_{k\ge1}, \quad \lambda_{k\nu j}=\rho_{k\nu j}^2,
$$
$$
\rho_{k\nu j} =\frac\pi{b_\nu-a_\nu}\Big(k-\frac{1-\delta_{N,\nu}-j\delta_{1,\nu}}{2}\Big) + o(1), \quad k\to\infty.
$$

Now let us study the asymptotical behavior of the functions $C(x,\lambda)$ and  $\Phi(x,\lambda).$ For our purposes it is sufficient to
consider $\rho\in\Omega_{\de}:=\{z:\arg z \in[\de,\pi-\de]\}$ and $x\in(0,b_1).$ For evaluating $C(x,\lambda)$ one can use relations
(\ref{24}) and (\ref{25}). However, on $(0,b_1)$ the function $C(x,\lambda)$ coincides with the classical cos-type solution (see, e.g.,
\cite{[5]}) and has the asymptotics
\begin{equation}  \label{16S}
C^{(\nu)}(x, \lambda) = \frac{(-i\rho)^\nu}2\exp(-i \rho x)[1], \quad x\in(0,b_1], \quad \rho\in\Omega_\delta, \quad \nu=0,1.
\end{equation}

Further, expanding $\Phi(x,\lambda)$ with respect to the system  $\{Y_1(x,\rho), Y_2(x,\rho)\},$ we get
\begin{equation}\label{14}
\Phi(x,\lambda)=A_{2k-1}(\rho) Y_1(x,\rho)+A_{2k}(\rho) Y_2(x,\rho), \quad x \in[a_k,b_k],\quad k=\overline{1,N},
\end{equation}
where the vector $A=\big(A_k(\rho)\big)_{k=1}^{2N}$ can be found from jump conditions (\ref{6}) and boundary conditions (\ref{7.1}).
Indeed, substituting (\ref{14}) into (\ref{6}) and (\ref{7.1}) we obtain $D_1A=(-1,0,\ldots,0)^T.$ Solving this system by Cramer's
formulae, we get
\begin{equation}\label{A1A2}
A_1(\rho)=\frac{q_{22} {\cal D}_0^{N-1} - q_{21} {\cal D}_1^{N-1}}{{\cal D}_1^N}, \quad A_2(\rho)=\frac{p_{21} {\cal D}_1^{N-1} - p_{22}
{\cal D}_0^{N-1}}{{\cal D}_1^N}.
\end{equation}

According to (\ref{27}), in $\Omega_\delta$ we have
$$
{\cal D}_j^k =(-1)^k(i\rho)^{k+j-1}\prod_{l=1}^{k-1} \alpha_{22}^{N-l}(\lambda)\exp(-i\rho\gamma_k)[1], \quad k=\overline{1,N}, \quad
j=0,1,
$$
which along with (\ref{alpha}), (\ref{pq}) and (\ref{A1A2}) give
$$
A_1(\rho)=\frac1{i\rho}[1], \quad A_2(\rho)=\frac1{i\rho}\exp(2i\rho b_1)[1].
$$
Substituting this into (\ref{14}) for $k=1$ and using (\ref{12}), we arrive at
$$
\Phi^{(\nu)}(x,\lambda) =(i\rho)^{\nu-1}\exp(i\rho x)\Big([1] +(-1)^\nu\exp(2i\rho(b_1-x))[1]\Big), \quad x\in[0,b_1],
$$
and, hence,
\begin{equation}\label{Phi}
\Phi^{(\nu)}(x,\lambda) =(i\rho)^{\nu-1}\exp(i\rho x)[1], \quad x\in[0,b_1).
\end{equation}
\\

\noindent {\bf 4. The inverse problem}
\\

Consider the following inverse problem.

\medskip
{\bf Inverse Problem 1.} Given $M(\lambda),$ find $q$ on $T.$

\medskip Using the ideas of the method of spectral mappings \cite{[6]} we prove the uniqueness theorem for the solution of Inverse
Problem~1. For this purpose together with the boundary value problem $L_1$ we consider a problem $\tilde L_1$ of the same form but with
another potential $\tilde q.$ We agree that if a certain symbol $\theta$ denotes an object related to $L_1,$ then this symbol with tilde
$\tilde\theta$ will denote the analogous object related to $\tilde L_1.$

\medskip
{\bf Theorem 1. }{\it If $M(\lambda)=\tilde M(\lambda),$ then $q=\tilde q$ on $T.$ Thus, specification of the Weyl function $M(\lambda)$
uniquely determines the potential $q.$}

\medskip
{\it Proof.} (i) At first, let us prove that $q$ and $\tilde q$ coincide on the segment $[0,b_1].$ For $x\in(0,b_1)$ consider the functions
$$
P_1(x,\lambda)=\tilde\Phi'(x,\lambda)C(x,\lambda)-\Phi(x,\lambda)\tilde C'(x,\lambda),\quad P_2(x,\lambda)=\Phi(x,\lambda)\tilde
C(x,\lambda)-\tilde\Phi(x,\lambda)C(x,\lambda).
$$
By virtue of the relation $C(x,\lambda)\Phi'(x,\lambda)-C'(x,\lambda)\Phi(x,\lambda)\equiv1,$ we have
\begin{equation}\label{P}
P_1(x,\lambda)\tilde C(x,\lambda)+P_2(x,\lambda)\tilde C'(x,\lambda)=C(x,\lambda).
\end{equation}
Further, from (\ref{16S}) and (\ref{Phi}) it follows that for each fixed $x\in(0,b_1)$
\begin{equation}\label{18}
P_1(x,\lambda)=1+O\Big(\frac1\rho\Big),\quad P_2(x,\lambda)=O\Big(\frac1{\rho^2}\Big),\quad \rho\to\infty,\quad \rho\in\Omega_{\delta}.
\end{equation}
On the other hand, using (\ref{9}) and the assumption of the theorem, we get
$$
P_1(x,\lambda)=C(x,\lambda)\tilde S'(x,\lambda)-\tilde C'(x,\lambda)S(x,\lambda),\; P_2(x,\lambda)=\tilde
C(x,\lambda)S(x,\lambda)-C(x,\lambda)\tilde S(x,\lambda),
$$
and consequently, for each fixed $x\in(0,b_1),$ the functions $P_1(x,\lambda)$ and $P_2(x,\lambda)$ are entire in $\lambda$ of order $1/2.$
By the Phragmen--Lindel\"of theorem and Liouville's theorem, asymptotics (\ref{18}) imply $P_1(x,\lambda)\equiv 1$ and
$P_2(x,\lambda)\equiv 0,$ which along with (\ref{P}) gives $C(x,\lambda)=\tilde C(x,\lambda)$ for $x\in(0,b_1)$ and, by continuity, for
$x\in[0,b_1].$ Then $q(x)=\tilde q(x)$ for $x\in[0,b_1].$

(ii) Let us prove by induction that the Weyl function uniquely determines the potential $q(x)$ on the other segments $[a_k,b_k],\;
k=\overline{2,N}.$ For this purpose we split $T$ into the union of the sets
$$
T_m:=\bigcup_{k=m}^N\,[a_k, b_k],\quad T_{m,0}:=\bigcup_{k=1}^{m-1}\,[a_k, b_k]
$$
and determine
\begin{equation}\label{20}
\Phi_m(x,\lambda):=\frac{\Phi(x,\lambda)}{\Phi'(a_m,\lambda)},\quad M_m(\lambda):=\Phi_m(a_m,\lambda),\quad m=\overline{2,N}.
\end{equation}
Then
\begin{equation}
M_m(\lambda)=\frac{\Phi(a_m,\lambda)}{\Phi'(a_m,\lambda)}.                  \label{21}
\end{equation}
According to (\ref{8}) and (\ref{20}), we have $\Phi_m'(a_m,\lambda)=1,\;\Phi_m(b_N,\lambda)=0,$ and consequently, the function
$M_m(\lambda)$ is the Weyl function for equation (\ref{4}) on $T_m.$

Fix $m=\overline{2,N},$ and suppose that we already proved that $q(x)\equiv \tilde q(x)$ for $x\in T_{m,0}.$ Then
$$
C(x,\lambda)\equiv\tilde C(x,\lambda),\quad S(x,\lambda)\equiv\tilde S(x,\lambda),\quad x\in T_{m,0}\,.
$$
Using (\ref{9}) along with the assumption of the theorem we get
\begin{equation}
\Phi(x,\lambda)\equiv\tilde\Phi(x,\lambda),\quad x\in T_{m,0}.           \label{22}
\end{equation}
Taking jump conditions (\ref{6}) for the function $y=\Phi(x,\lambda)$ and (\ref{7}), (\ref{22}) into account we infer
$$
\Phi(a_m,\lambda)\equiv\tilde\Phi(a_m,\lambda), \quad \Phi'(a_m,\lambda)\equiv\tilde\Phi'(a_m,\lambda).
$$
Together with (\ref{21}) this yields
$$
M_m(\lambda)=\tilde M_m(\lambda).                                     
$$
Repeating the arguments from the first part of the proof, we conclude that $q(x)\equiv\tilde q(x)$ for $x\in[a_m,b_m],$ which finishes the
proof. $\hfill\Box$

\medskip
Using Theorem~1 we also prove the uniqueness theorem for the inverse problem of recovering the potential from the spectra of two boundary
value problems $L_0$ and $L_1.$

\medskip
{\bf Inverse problem 2. } Given $\{\lambda_{nj}\}_{n\ge 1},\, j=0,1,$ find $q$ on $T.$

\medskip
Hadamard's factorization theorem gives
$$
\Delta_j(\lambda) = C_j p_j(\lambda), \quad p_j(\lambda) = \lambda^{s_j}\prod_{\lambda_{nj}\ne0} \Big(1 -
\frac{\lambda}{\lambda_{nj}}\Big), \quad j=0,1,
$$
where $C_j$ is a non-zero complex constant, while $s_j$ is the multiplicity of the zero eigenvalue in the spectrum $\{\lambda_{nj}\}_{n\ge
1}.$ Consider the functions
$$
g_j(\lambda) = (-i\rho)^{2N+j-3} \prod_{l=1}^{N-1}(a_{l+1}-b_l)^2 f_j^N(\rho)f_1^1(\rho)\prod^{N-1}_{l=2}f_0^l(\rho), \quad j=0,1.
$$
By virtue of \eqref{Delta}, the following limits exist:
$$
\lim_{\lambda\to-\infty}\frac{\Delta_j(\lambda)}{g_j(\lambda)}=1, \quad j=0,1,
$$
and, hence,
$$
C_j = \lim_{\lambda\to-\infty}\frac{g_j(\lambda)}{p_j(\lambda)}.
$$
Thus, the characteristic functions $\Delta_j(\lambda)$ are uniquely determined by their zeros $\{\lambda_{nj}\}_{n\ge 1}.$ Using formula
(\ref{10}) along with Theorem~1 we arrive at the following uniqueness theorem for Inverse problem~2.

\medskip
{\bf Theorem 2. }{\it Specification of the spectra $\{\lambda_{nj}\}_{n\ge 1},$ $j=0,1,$ uniquely determines the potential $q(x)$ on $T.$}
\\

{\bf Acknowledgment.} This work was supported in part by Grant 1.1660.2017/4.6 of the Ministry of Education and Science of the Russian
Federation and by Grant 19-01-00102 of Russian Foundation for Basic Research.\\

\end{document}